\documentclass[reqno,twoside,11pt]{amsart}

\usepackage{amsmath,amsfonts,calrsfs,fullpage,amssymb,color,verbatim,eucal,yfonts,mathrsfs}
\usepackage[bookmarksnumbered, colorlinks, plainpages]{hyperref}
\usepackage{amsfonts,latexsym,amstext}
\usepackage{amsmath, amssymb,xcolor}
\usepackage[normalem]{ulem}

\newtheorem{theorem}{Theorem}[section]

\newtheorem{proposition}[theorem]{Proposition}
\newtheorem{definition}[theorem]{Definition}
\newtheorem{example}[theorem]{Example}

\newtheorem{Assumption}[theorem]{Assumption}
\newtheorem{remark}[theorem]{Remark}
\newtheorem{corollary}[theorem]{Corollary}

\numberwithin{equation}{section}

\def\sqr#1#2{{\vcenter{\vbox{\hrule height .#2pt \hbox{\vrule
 width .#2pt height#1pt \kern#1pt \vrule
width .#2pt} \hrule height .#2pt}}}}
\def\square{\mathchoice\sqr54\sqr54\sqr{4.1}3\sqr{3.5}3}

\def\ds{\begin{displaystyle}}
\def\eds{\end{displaystyle}}

\def\<{\langle }
\def\>{\rangle }

\newcommand{\e}{\epsilon}

\newcommand{\E}{\mathbb{E}}
\newcommand{\R}{\mathbb{R}}
\newcommand{\F}{\mathcal{F}}
\def\Dim{\noindent\hbox{{\bf Proof.}$\;\; $}}
\def\finedim{{\hfill\hbox{\enspace${ \square}$}} \smallskip}

\allowdisplaybreaks

\begin{document}

\title{On approximations of stochastic optimal control problems with an application to climate equations}

\author[Franco]{Franco Flandoli \textsuperscript{1}}

\address{\textsuperscript{1}Scuola Normale Superiore, Piazza dei Cavalieri, 7, 56126 Pisa, Italia}

\email{franco.flandoli@sns.it}

\author[Biba]{Giuseppina Guatteri \textsuperscript{2}}

\address{\textsuperscript{2}Dipartimento di Matematica,
Politecnico di Milano, Piazza Leonardo da Vinci 32,
20133 Milano, Italia.}

\email{giuseppina.guatteri@polimi.it}

\author[Umberto]{Umberto Pappalettera \textsuperscript{3}}

\address{\textsuperscript{3}Fakultat f\"ur Mathematik, Universit\"at Bielefeld, D-33501 Bielefeld, Germany }

\email{upappale@math.uni-bielefeld.de}

\author[Mario]{Gianmario Tessitore \textsuperscript{4}}

\address{\textsuperscript{4} Università di Milano-Bicocca, Dipartimento di Matatematica e Applicazioni, via Roberto Cozzi 55,
20125 Milano, Italia.}

\email{gianmario.tessitore@unimib.it}


\begin{abstract}
The paper is devoted to the optimal control of a system with two time-scales, in a regime when the limit equation is not of averaging type but, in the spirit of Wong-Zakai principle, it is a stochastic differential equation for the slow variable, with noise emerging from the fast one. It proves that it is possible to control the slow variable by acting only on the fast scales. The concrete problem, of interest for climate research, is embedded into an abstract framework in Hilbert spaces, with a stochastic process driven by an approximation of a given noise. The principle established here is that convergence of the uncontrolled problem is sufficient for convergence of both the optimal costs and the optimal controls. This target is reached using Girsanov transform and the representation of the optimal cost and the optimal controls using a Forward Backward System. A challenge in this program is represented by the generality considered here of unbounded control actions.
\end{abstract}

\maketitle

\section{Introduction}
In this paper we are concerned with optimal control problems associated to stochastic equations in abstract Hilbert spaces \cite{DPZ}, and their convergence.
More precisely, we introduce a family of stochastic equations indexed by a parameter $\e \in (0,1)$, driven by a stochastic process obtained as \textit{approximation} of some given noise $W$ via some general approximation map $\Gamma^\e$.
Then, we solve a control problem for every $\e \in (0,1)$ and we are interested in understanding if the convergence of both the optimal costs and the optimal controls holds true as $\e \to 0$.

In order to answer this question, we develop a general framework for studying approximations of stochastic optimal control problems.

A part from very natural technical assumptions, the only hypothesis on the approximation maps $\Gamma^\e$ is the validity of some Wong-Zakai type of convergence. That is, we assume that the solution of the uncontrolled equation driven by the approximation of the noise $\Gamma^\e(W)$ converges in probability, as $\e \to 0$, towards the solution of the uncontrolled equation driven by $W$ (cf. Assumption \ref{ipotesi limite forward}).
The key result of this work is that convergence of the uncontrolled problems is sufficient for convergence of both the optimal costs and the optimal controls (cf. Theorem \ref{main_genenral})

{In this paper the above goal is  achieved through the representation of  the optimal cost and the optimal controls using a Forward Backward System  of Stochastic Differential Equation (FBSDEs) (see e.g. the system \eqref{stato eps new} here).
This technique has been widely used in the last twenty-five years both in finite and in infinite dimensional framework (see for instance \cite{YongZhou} or \cite[Chapter 6]{FaGoSw} and references within). It has the advantage to characterize not only the optimal state of a stochastic control problem but also the optimal feedback law, without requiring regularity of the value function. This is, in extreme synthesis the reason why we are able to obtain our main abstract convergence result, see Theorem \ref{main_genenral} and, in particular, the convergence of optimal controls stated in it.}


{ In this work, we consider the case in which the running cost is quadratic and coercive with respect to the control variable, while exhibiting bounded behavior in the state variable (cf. Assumption \ref{hyp-control}). This choice of allowing ``unbounded'' control actions introduces significant technical challenges in the development of the Forward-Backward Stochastic Differential Equation (FBSDE) approach to optimal control problems. First of all it interacts with the necessity of adopting a weak formulation for the control problem. Indeed the final convergence argument (cf. Section \ref{sec form contr}) works if uncontrolled state equations  refer to the same stochastic framework. Consequently, we are led to express the control problem in a weak form. This implies that a rigorous formulation of the problem, along with the characterization of the class of admissible controls, involves a change of probability that necessitates the introduction of a localization argument. (cf. Definition \ref{Controlli Ammissibili eps} and Proposition \ref{well posedness cost}).}

{ In addition, the Hamiltonian non-linearity \(\psi\), introduced in Section \ref{sec hamiltonian} , which drives the backward equation in the FBSDE system (cf. Equation \eqref{stato eps new}), is non-Lipschitz with respect to its second variable, denoted by $Z$. To address this point, we adapt the techniques developed in  \cite{Kobilanski} and \cite{BriConf}, taking profit, in particular, of the specific properties of BMO martingales (cf. Section  \ref{sec BMO} here  and \cite{Kazamaki}). The appropriate use of this class of martingales, along with the corresponding estimates, constitutes a crucial element in the proof of our main general result, c.f. Theorem \ref{main_genenral}.
 In synthesis,  the combination of a coercive quadratic cost function and a non-Lipschitz Hamiltonian introduces considerable complexities. These necessitate the use of advanced stochastic analysis techniques, notably those concerning BMO martingales, to ensure the well-posedness of the optimal control problems  and their convergence, within the weak formulation framework .}

Our main motivation for studying general approximations of stochastic optimal control problems comes from the desire of understanding the behaviour of controlled slow-fast systems of stochastic equations $(X^\e,Q^\e)$, depending on a small parameter $0<\e\ll 1$.
Indeed, in certain prototypical situations, the slow component $X^\e$ of the system converges as $\e \to 0$ towards a limiting closed equation. Here the term ``closed'' refers to the fact that the equation for the limit $\hat{X}$ no longer depends on the fast variable. 
In this case, we intend to study a control problem for the systems $(X^\e,Q^\e)$ and for the limit equation $\hat{X}$.
A natural question is whether the control problem for $(X^\e,Q^\e)$ can be solved at every $\e >0$, and whether or not the solutions of the control problems converge as $\e \to 0$ to a solution of the control problem for $\hat{X}$. 
The relevance of this problem becomes clear in view of the interpretation of slow-fast systems as general models of climate-weather interaction (see next subsection for additional details). 
With the lens of this interpretation, the convergence of control problems translate into the following question: Is it possible to "control" the evolution of the climate by acting only at meteorologic scales?  

We believe this setting is robust enough to be amenable to further generalizations of the control problems \eqref{statopapcontr} and \eqref{statopapreducedcontr}, cf. the discussion in \autoref{sec:Examples}.

\subsection{A Motivating Example: Climatic Model}\label{secFlandoliPappalettera}
Let us start with a motivating example. 
We consider a slow-fast system having the following form: 
 \begin{equation} \label{statopap} \left \{
     \begin{array}{ll}
        d X^\e_t = A X^\e_t dt + {b}(X^\e_t) dt  + \sigma( X^\e_t) Q^\e_t dt,   &  t \in [0,T], \\ 
         X^\e(0)= x_0, & \\
         d Q^\e_t= -\dfrac{1}{\e}  Q^\e_t dt + \dfrac{1}{\e} Gd{W}_t, &  t \in [0,T], \\ 
         Q^\e(0)=0. &
     \end{array} \right.
 \end{equation}
 Solutions of \eqref{statopap} are pairs of stochastic processes $(X^\e$, $Q^\e)$, where the ``slow'' component $(X^\e)$ takes values in a Hilbert space $K$ and  the ``fast'' component $(Q^\e)$ takes values in a Hilbert space $H$.
 We denote $|\cdot|_K$ and $|\cdot|_H$ the norms on these spaces, and $\langle \cdot ,\cdot \rangle_K$ and $\langle \cdot ,\cdot \rangle_K$ the inner products.
 For simplicity we assume $x_0 \in K$ given and deterministic.
 In the lines above, $(W_t)$ is a cylindrical Wiener process on a complete probability space $(\Omega, \mathcal{F},\mathbb{P})$ with complete and right-continuous filtration $(\mathcal{F}^W_t)$, and $G$ is a Hilbert-Schmidt operator on the Hilbert space $H$ with $Gv=\sum_{i=1}^\infty \lambda_i \langle e_i,v\rangle_H e_i$ for every $v\in H$, where $\sum_{i=1}^\infty \lambda_i^2 <\infty$ and $(e_i)_{i\in\mathbb{N}}$ is a orthonormal basis in $H$.
 Finally, $A:K \rightarrow K$ is linear continuous, $b:K\rightarrow K$ is Lipschitz and  $\sigma:K \rightarrow \mathcal{L}(H; K)$ is of class $C^2_b$. 
 The maps $b$ and $\sigma$ can be expressed in terms of their coordinates $b(x)_j := \langle b(x) , f_j \rangle_K$ and $(\sigma(x)e_m)_j:=\sigma^{j,m}(x) :=\langle\sigma(x)e_m, f_j\rangle_K$, where $(f_j)_{j\in\mathbb{N}}$ is a orthonormal basis in $K$ and $j,m\in \mathbb{N}$.
 
This kind of slow-fast systems have been extensively studied in pure and applied mathematics.
Among other important question, one is naturally led to ask what the behaviour of this system is in the limit of infinite separation of scales $\e \to 0$. 
Heuristically, the fast oscillations of the process $Q^\e$ prevent it from converging as a genuine function, and convergence of $Q^\e$ usually holds in a space of distributions with respect to time. 
On the other hand, the slow component $X^\e$ can converge as a function but its limiting dynamics should retain information about the statistics of $Q^\e$.
When the limit $\hat{X} := \lim_{\e \to 0} X^\e$ solves a closed equation, we say that the limiting equation for $\hat{X}$ is a \textit{stochastic model reduction} of \eqref{statopap}.

The first rigorous examples of stochastic model reduction of finite dimensional equations are due to Kurtz \cite{kurtz} and Majda, Timofeyev, and Vanden Eijnden \cite{MTV}.
In particular, the latter successfully gave a stochastic model reduction of the truncated Barotropic Equations, identifying the slow variable $X^\e$ as a quantity evolving on \textit{climatic} time-scale and the slow variable $Q^\e$ as a quantity evolving on \textit{meteorologic} time-scale. 
The small constant $\e>0$ represents the ratio between the speed of the evolution at these different time-scales.
A similar interpretation was given in \cite{AsFlaPap}.

Under the previous assumptions on \eqref{statopap} and assuming $K$ finite dimensional, in \cite{AsFlaPap} it is proved that, as $\e$ goes to zero, the sequence $(X^\e)$ converges in probability in the $C([0,T], K)$ norm towards the  solution $(\hat{X})$ of the ``reduced'' equation 
 \begin{equation} \label{statopapreduced}
 \left \{
     \begin{array}{ll}
        d\hat{X}_t = A\hat{X}_t dt + \hat{b}(\hat{X}_t) dt  + \sigma( \hat{X}_t)GdW_t,  &  t \in [0,T], \\ 
         \hat{X}(0)= x_0,
     \end{array} \right.
 \end{equation}
 where 
 \begin{align*}
    (\hat{b}(x))_i:=(b(x))_i+\dfrac{1}{2}\sum_{m=1}^{\infty} \lambda_m^2\sum_{j=1}^{d} D_j\sigma^{i,m}(x)\sigma^{j,m}(x), \quad x \in K, i \in \mathbb{N}. 
 \end{align*}
 Notice that under the present assumptions $\hat{b}:K \rightarrow K$ is Lipschitz.

 {It is wort noticing that, differently to  ``standard'' two-scales stochastic models where the fast evolution equation is obtained by a simple change of the  time-scale with  ratio $\epsilon$ (for the controlled version of such systems see, e.g.  \cite{KabPer}, \cite{AlvBar} \cite{guatess1}, \cite{guatess2}, and \cite{Swiech2021}) here the oscillations induced by the noise in the fast equation is magnified by a factor $1/\sqrt{\epsilon}$. As a matter of fact the two class of slow-fast models show a very different behaviour in the limit. Here a new noise term and a correction drift appears in the reduced equation while in the other case the reduced equation is obtained by ``averaging''  the original coefficients with respect to a suitable ``invariant measure''}

It should also be pointed out that equations of the form \eqref{statopapreduced} have already appeared in the study of climate since the seminal work of Hasselmann \cite{Ha76} on \textit{stochastic climate models}.
Indeed, Hasselmann proposes a general stochastic model to predict the evolution of quantity on \textit{climatic} time-scales, without referring to any particular specification of the coefficients $A$, $\hat{b}$, $\sigma$ of \eqref{statopapreduced}.
The deep aspect of Hasselmann proposal is that a (small intensity) noise should be taken into account for a more correct description of the system. 

In a second moment, the general theory of stochastic climate models has been  specialized to particular systems, possibly adding \textit{ad hoc} assumptions on the coefficients. 
To mention a few works in this direction, let us cite \cite{Ha77} on sea-surface temperature anomalies and thermocline variability, \cite{DelSartoFranco} on a energy balance model addressing temperature fluctuations due to rising carbon dioxide levels, \cite{ButoriLuongo} on magneto-hydrodynamics models, \cite{FPT} on random attractors, and \cite{Lucarini} on climatic tipping points. 

\subsection{Controlled climatic model}
We wish to study a controlled version of this model, with control acting at the meteorologic scale. 
Namely, fixed {an Hilbert space $U$} and given a progressively measurable control process $u$ taking values in $U$, we consider the system
 \begin{equation} \label{statopapcontr} \left \{
     \begin{array}{ll}
        d X^{\e,u}_t = AX^{\e,u}_t dt + {b}(X^{\e,u}_t) dt  + \sigma( X^{\e,u}_t) Q^{\e,u}_t dt,   &  t \in [0,T], \\ 
         X(0)= x_0, & \\
         d Q^{\e,u}_t= -\dfrac{1}{\e}  Q^{\e,u}_t dt +\dfrac{1}{\e} Gr(u_t)dt+  \dfrac{1}{\e} Gd{W}_t, &  t \in [0,T], \\ 
         Q(0)=0, &
     \end{array} \right.
 \end{equation}
 where $r: U\rightarrow H$ is a Lipschitz map.
 We also introduce a controlled reduced equation, namely
 \begin{equation} \label{statopapreducedcontr} \left \{
     \begin{array}{ll}
        d\hat{X}^u_t = A\hat{X}_t dt + \hat{b}(\hat{X}_t) dt + \sigma( \hat{X}^u_t)Gr(u_t)dt  + \sigma( \hat{X}^u_t)GdW_t,  &  t \in [0,T], \\ 
         \hat{X}(0)= x_0.
     \end{array} \right.
 \end{equation}

The control problems above come with the two cost functionals: $J^\e$, related to system \eqref{statopapcontr} and $J$, related to equation \eqref{statopapreducedcontr}, that  we assume both quadratic and coercive in $u$.
In details 
  the costs are given by
 $$J^{\epsilon}(x_0,u) := \mathbb{E}\left[\int_0^T l(X^{\epsilon,u}_s,u_s)ds +h(X^{\epsilon,u}_T)\right] \; \hbox{ and } 
 J(x_0,u) := \mathbb{E}\left[\int_0^T l(X^{u}_s,u_s)ds +h(X^{u}_T)\right]. $$
so that the functional above are well defined.

 For the precise assumptions on $l$ and $h$, as well as for the definition of the class of admissible controls, we refer to Assumption \ref{hyp-control}, Definition \ref{Controlli Ammissibili eps}, and Theorem \ref{teo esist limite} below.
 
 Our main result goes as follows, see also Theorem \ref{main_genenral} for a precise statement.
 We prove that:

$i$) The control problems \eqref{statopapcontr} and \eqref{statopapreducedcontr} admit an optimal control, denoted respectively $\underline{u}^\e$ and $\hat{u}$;

$ii$) The optimal controls are square integrable;

$iii$) As $\e \to 0$, the optimal costs converge: $J^\e(x_0,\underline{u}^\e) \to J(x_0,\hat{u})$;

$iv$) As $\e \to 0$, the optimal controls converge:
$ \displaystyle \E \int_0^T | \underline{u}^\e _t- \hat{u}_t| ^2 \, dt \to 0$



{It is perhaps worth noticing that although we start by approximating problems with control acting at meteorological time scale we end up with a limit reduced problem with control acting at climatic time scale}

{We hope that  this work, devoted to a simplified model,  may serve as a useful starting point for examining the behaviour of optimal controls in related, more realistic, contexts.}

\section{A general framework for approximation of stochastic optimal control problems}

In this section we state the control problem we are going to study.
We will introduce a  weak formulation of the  problem that will be particularly suitable for our  purposes. Indeed, it allows to formulate all control problems on the same stochastic basis; this fact, together with the representation of the optimal cost and optimal control by the solution of a backward stochastic differential equation (BSDE), will allow to show the convergence, as $\e \to 0$, of both the optimal costs $J^\e(x_0,\underline{u}^\e) \to J(x_0,\hat{u})$ and the optimal controls $\underline{u}^\e \to \hat{u}$.
We underline  that, being the controls unbounded, the weak formulation of the control problem can not be obtained by a standard application of Girsanov transform, and we have to proceed by localization.

\subsection{Settings} \label{settings}
{
Let us reprise the notation of the previous motivating examples, although precise working assumptions in the more general framework will be stated later.
We fix  $H$, $K$ and $U$ real separable Hilbert spaces (in our motivating model $H$  hosts the fast variables, $K$ should host the slow variables of the system, while $U$ will describe the actions of the control processes). Also recall that $(W_t)_{t \geq 0}$ is a $H$-valued Wiener process on a complete probability space $(\Omega,\mathcal{F},\mathbb{P})$, with complete and right-continuous filtration $(\mathcal{F}^W_t)_{t \geq 0}$.}

{
Given an arbitrary Banach space $E$ and $p \in [1,\infty)$ let us denote $L^{p,loc}_W(\Omega \times [0,T], {E})$ the space of $(\mathcal{F}^W_t)$-progressively measurable stochastic processes in 
\begin{align*}
    L^{p,loc}(\Omega \times [0,T]; E) 
    := 
    \left\{ \Phi : \Omega \times [0,T]\to E \,:\, \int_0^T | \Phi_t |_{E}^p dt < \infty \quad \mathbb{P}\mbox{-a.s.}\right\}.
\end{align*}}

{We will need also the spaces  $L^2 _W(\Omega \times [0,T]; E) $ of square integrable  progressive measurable processes $\Phi : \Omega \times [0,T]\to E$ verifying
\begin{align*} |\Phi|^2_{    L^{2}_W(\Omega \times [0,T]; E) }
    :=  \E \int_0^T | \Phi_t |_{E}^2 dt < \infty
\end{align*}
and, for $p\in [1,\infty]$ the spaces  $L^p_W(\Omega;C([0,T];E))$
       of
     progressive measurable processes $Y$ with continuous paths in $E$, such
    that the norm $|Y|_{L^p_W(\Omega;C([0,T];E))}:=
    |\sup _{s \in [0,T]} |Y_s|_E|_{L^p(\Omega)}$
    is finite, that
    the subspace of
    predictable processes $Y$ with continuous paths in $E$.}
 
{ We finally denote by $\mathcal{I}$ the space of $H$-valued continuous  It\^o-semimartingales of the form
 \begin{equation} \label{eq:I}
     I_t=\int_0^t \Phi_s ds+\int_0^t \Psi_s dW_s,
 \end{equation}
with $\Phi \in L^{1,loc}_W(\Omega \times [0,T]; H)$ and $ \Psi \in  L^{2,loc}_W(\Omega \times [0,T]; L_2(H))$ where $L_2(H)$ stands for the space of Hilbert-Schmidt operators  
 from $H$ to $H$. }
 
We introduce a class of functionals $(\Gamma^\e) _{\e >0}$ from $\mathcal{I}$ to the class of
càdlàg processes $\Omega\times[0,T]\rightarrow K$ we assume the following
\begin{Assumption}\label{ipotesi su Gamma}  $(\Gamma^\e)$ is  an $(\mathcal{F}^W_t)_{t\geq 0}$-adapted process  and its law  only depends on the law of $I$.
    \end{Assumption}
For the sake of the presentation, let us point out that one could think of the family of functionals $(\Gamma^\e) _{\e >0}$ as some \textit{adapted} approximation of the noise $I$, e.g. {adapted} piecewise linear interpolation, convolution, or coloured-in-time approximation à la Ornstein-Uhlenbeck.

\medskip


\subsection{State Equations}

The first class of state equations corresponds to a regularization of the noises induced by the functionals $(\Gamma^\e)$, namely:
{\renewcommand{\theequation}{S$_\e$}
 \begin{equation}\label{stato eps} \left \{
     \begin{array}{ll}
        d X^\e_t =  A X^\e_t dt +b(X^\e_t) dt + \sigma( X^\e_t)   \Gamma ^\e [G W] _t dt, &  t\in (0,T], \\ \\
         X^\e(0)= x_0. & 
     \end{array} \right.
 \end{equation}}
The second class corresponds to stochastic equations with white-in-time noises: 
 {\renewcommand{\theequation}{S}
 \begin{equation} \label{stato} \left \{
     \begin{array}{ll}
        d \hat{X}_t = ( A \hat{X}_t + \hat{b}(\hat{X}_t)) \, dt + \sigma( \hat{X}_t) G  dW_t, &  t\in (0,T], \\ \\
         X(0)= x_0. & 
     \end{array} \right.
 \end{equation}}
 \setcounter{equation}{0}

\begin{example}
The motivating example of the introduction can be rephrased within this general framework by definining
    \begin{equation} \label{eq:Gamma}
    \Gamma^\e[I]_t := \frac{1}{\e}\int_0^t e^{(t-s)A}dI_s=\frac{1}{\e}\int_0^t e^{(t-s)A}\Phi_sds+\frac{1}{\e}\int_0^t e^{(t-s)A}\Psi_sdW_s.
\end{equation}
Indeed, the cost functionals $J^\e$ and $J$ only involve the slow component of the solution to systems \eqref{statopapcontr} and \eqref{statopapreducedcontr}, so \eqref{statopapcontr} and \eqref{statopapreducedcontr} can be replaced by \eqref{stato eps} and \eqref{stato} without loss of useful information.
\end{example}
  
In both \eqref{stato eps} and \eqref{stato}, the coefficients satisfy the following: 
\begin{Assumption}\label{ass-state} We assume the following:
\begin{enumerate}
\item[(Hp 1-A)] $A:D(A)\to K $ is a  (possibly unbounded) linear operator with domain $D(A) \subseteq K$ that generates a $C_0$- semigroup $(e^{tA})_{t\geq 0}$.
    \item[(Hp 1-b)]  $b$ and $\hat{b}$ are Lipschitz maps from $K$ to $K$ and we fix a constant 
    $L_b >0$ such that
    \begin{align}
        \label{def b}
   |b(x)-b(y)|_K &\leq  L_b \, |x-y|_K, \qquad\qquad\qquad \ \   \forall x,y \in K;
\end{align}
and the same holds for $\hat{b}$.
 \item[(Hp 1-$\sigma$)] $\sigma$ is a Lipschitz map from $K$ to $L(H;K)$ 
 and we fix a constant 
    $L_\sigma >0$ such that
    \begin{align}
        \label{def sigma}
   |\sigma(x)-\sigma(y)| _{L(H;K)} &\leq  L_\sigma |x-y| _K, \qquad\qquad\qquad \ \ \forall x,y \in K;
\end{align}
\end{enumerate}
\end{Assumption}

The following existence and uniqueness result is a consequence of straightforward fixed point arguments.
\begin{theorem}\label{esistenza stato}
    Under Assumption \ref{ass-state}, for every $\e >0$ there exists an adapted process $X^\e$ with continuous trajectories solving equation \eqref{stato eps} in a mild-pathwise sense; that is, such that for $\mathbb{P}$ almost every $\omega\in \Omega$, it holds for all $t\in [0,T]$
$$ X^{\e}_t(\omega)
=
x_0+\int_0^t e^{(t-s)A}\Big(  b(X^{\e}_s(\omega))+\sigma(X^{\e}_s(\omega))\Gamma^\e[GW]_s(\omega) \Big) ds.$$
Moreover, there exists a unique mild solution $\hat X$ of equation \eqref{stato} that belongs to $L^p_W(\Omega; C([0,T];K))$ for all $p>1$.

\end{theorem}


\medskip


\subsection{Controlled equations}
Next, let us introduce the controlled equations we are going to study in this general framework. 
Let $X^{\e,u}$ solve
\begin{equation}\label{stato eps u} \left \{
     \begin{array}{ll}
       d X^{\e,u}_t =  A X^{\e,u}_t dt +b(X^{\e,u}_t) dt + \sigma( X^{\e,u}_t)  (\Gamma ^\e [G (W+\int_0^{\cdot}r(X^{\e,u}_s,u_s)ds)])_t dt, &  t \in (0,T], \\ \\
         X^\e(0)= x_0, & 
     \end{array} \right.
 \end{equation}
and let $X^{u}$ solve
 {\renewcommand{\theequation}{S}
 \begin{equation} \label{stato u} \left \{
     \begin{array}{ll}
        d X^u_t = ( A X^u_t + \hat{b}(X^u_t)) \, dt + \sigma( X^u_t) G  \, dW_t + \sigma( X^u_t) G  \, r(X^u_t,u_t)dt, &  t\in (0,T], \\ \\
         X(0)= x_0. & 
     \end{array} \right.
 \end{equation}}

Formally speaking, our purpose is to minimize the cost functionals (formally written) over all the admissible controls $u$
\begin{equation}\label{def funzionale controllo}
  J^{\epsilon}(x_0,u)=\mathbb{E}\left[\int_0^T l(X^{\epsilon,u}_s,u_s)ds +h(X^{\epsilon,u}_T)\right] \; 
\hbox{ and } 
J(x_0,u)=\mathbb{E}\left[\int_0^T l(X^{u}_s,u_s)ds +h(X^{u}_T)\right].     
\end{equation}
However, a precise formalization of the control problem will be given later.
For the time being, let us state our main assumptions on the functions $r$, $l$ and $h$.
\begin{Assumption}\label{hyp-control} We assume that:
\begin{enumerate} 
\item[(Hp 2-r)]  $r:K\times U \to H$ measurable such that for some constants $M_r, L_r >0$
    \begin{align}
        \label{def r}
   |r(x,u)|_H  &\leq  M_r (1+ |u|_U),  \qquad\qquad\qquad \  \ \forall x \in K, u \in U;
   \\
        \label{def rr}
   |r(x,u)-r(y,u)|_H &\leq  L_r( |x-y|_K\wedge 1)(|u|_U+1), \qquad\qquad\qquad \  \  \, \forall x,y \in K, u \in U.
\end{align}
Moreover, we assume that there exists $u_\star\in U$ such that $r(x,u_\star)=0$ for all $x\in H$. 
\item[(Hp 2-l)] $ l: K\times U \to  \R $ is a measurable map such that for some constants $M_l, {m}_l , c_l>0$   
\begin{align}
    \label{def l}
 {m}_l |u|_U ^2 - c_l \leq  l(x,u) & \leq  M_l(1+|u|_U^2), \qquad\qquad  \forall x \in K, u \in U;
 \\
    |l(x,u)-l(y,u)| &\leq  L_l |x-y| _K, \qquad\qquad \, \forall x,y \in K, u \in U.
\end{align}
Notice that the above implies that, for a suitable constant $C_l>0$
\begin{equation}\label{def g}
 |l(x,u)| \leq  C_l (1+|u|_U^2)  \qquad\qquad  \forall x \in K, u \in U
\end{equation}
 and hence we deduce that there exists a constant $C >0$ such that
\begin{equation}\label{def ll}
 |l(x,u) -l(y,u)| \leq  C ( |x-y|_K\wedge 1) (1+|u|^2_U)  \qquad\qquad  \forall x,y \in K, u \in U
\end{equation}

\item[(Hp 2-h)] $ h:  K \to  \R $ such that for some constant $M_h >0$  
\begin{align}
\label{def h}
   |h(x)|  &\leq  M_h, \qquad\qquad\quad \forall x \in K; \\
   |h(x)-h(y)| &\leq  L_h |x-y|_K,  \qquad\qquad\,  \forall x,y \in K.
\end{align}
\end{enumerate}

\end{Assumption}
\begin{remark} Notice that if $r(x,u)=r_0(x)u$ 	with $r_0$ bounded and Lipschitz then assumption {\em (Hp 2-r)} holds with $u_\star=0$.
\end{remark}

\subsection{Rigorous formalization of the control problem}\label{sec form contr}
We start by considering, for $\e>0$, the formal cost functional 
\begin{align*}
    J^{\epsilon}(x_0,u)
    =
    \mathbb{E}\left[\int_0^T l(X^{\epsilon,u}_s,u_s)ds +h(X^{\epsilon,u}_T)\right].
\end{align*}
If we assume the following boundedness condition on the controls
\begin{align*}
    \int_0^T|u_t|_U^2 dt \leq c < \infty, \quad  \mathbb{P} \mbox{ almost surely},
\end{align*} 
then a straightforward application of Girsanov transform, together with the fact that the law of the solution to equations \eqref{stato eps} does not depend on the specific stochastic basis, yields: 
$$
    J^\e(x_0 , u)= \E \Big[ \mathcal{E}_T ( r(X^\e,u)) \Big(\int_0^T l(X^\e_s,u _s)\, ds +  h(X^\e_T)\Big)\Big],
$$
    where  $\mathcal{E}_t (r(X^\e,u)):= \displaystyle\exp{\left\{ \int_0^ t   r(X^\e_s,u_s) \, d W_s -\frac{1}{2} \int_0^ t   |r(X^\e_s,u_s)|^2 \, ds \right\}}$ and, we recall, $(X^\e)$ solves the uncontrolled evolution eqaution \eqref{stato eps}.
    
However, while the requirement $u\in L^2_W(\Omega\times[0,T];U)$, is necessary in view of the quadratic behaviour with respect to $u$ of the running cost $l$ (see Assumption \eqref{def l} ,the requirement $\int_0^T|u_t|_U^2 dt \leq c$ seems artificially strong to be imposed in this context. 
In particular, optimal controls in the class $L^2_W(\Omega\times[0,T];U)$ do not have to satisfy it, see for instance the form of the optimal feedback $\hat{u}$ in the Example \ref{esempio-feedback}.


{Hence, we drop the assumption of $\mathbb{P}$-essential boundedness of  $L^2([0,T]; U)$ norm of the control trajectories. This, together with the choice of formulating our control problem in the weak  probabilistic form  causes several  technical difficulties.
Already the identification of  the  class of admissible controls is not trivial. 
Let us give the following definition that seems natural:  }
\begin{definition}\label{Controlli Ammissibili eps}
For every $u \in L^{2,loc}_W(\Omega\times [0,T]; U)$, we set the following.
\begin{enumerate}
    \item Let $\tau_n$, $n \in \mathbb{N}$, be a sequence of stopping times defined by
$$
\tau_n:= \inf{\left\{  t\geq 0: \int_0^t |u_s|_U^2 \, ds \geq n \right\}}.
$$
    \item Let $u_\star$ be such that $r(x,u_\star)= 0$ for every $x \in H$ (see (Hp 2-r) in Assumption \ref{hyp-control}) and let $u^n$ denote the control
$$
    u^n_s:= u_s I_{\{ 0 \leq s\leq \tau_n \wedge T \}} + u_\star I_{\{ \tau_n \wedge T < s \leq T \}},
    \quad s \in [0,T].
$$
\item Let $\mathcal{E}_t ( r(X^\e,u^n))$  be the exponential martingale 
$$
\mathcal{E}_t ( r(X^\e,u^n)) := \displaystyle\exp{\left\{ \int_0^ t   r(X^\e_s,u_s^n) \, d W_s -\frac{1}{2} \int_0^ t   |r(X^\e_s,u_s^n)|^2 \, ds \right\}}.
$$
\end{enumerate}

 We define the space of admissible controls $\mathcal{U}^\e_{ad}$ as the space:
 \begin{equation}\label{limite controlli def eps}
 \mathcal{U}^\e_{ad} 
 :=   
 \Big\{u \in L^{2}_W(\Omega\times [0,T]; U) 
 \,:\, 
 \sup_{n \in \mathbb{N}} \E \Big( \mathcal{E}_{T} ( r(X^\e,u^n)) \int_0^{T\wedge \tau^n} |u_s|_U^2 \,ds\Big) < \infty \Big\}.
 \end{equation}\end{definition}

Notice that $\mathcal{U}^\e_{ad} $, in the definition above, may in principle depend on $ x_0$.

\begin{remark}
    Notice that, in view of the fact that $r(x,u_\star)=0$ we have:
$$ \mathcal{E}_{T} ( r(X^\e,u^n))= \mathcal{E}_{T\wedge \tau^n} ( r(X^\e,u^n))=\mathcal{E}_{T\wedge \tau^n} ( r(X^\e,u)). $$
\end{remark}

\begin{remark} \label{rem lim}The sequence $\displaystyle n \mapsto \E \Big( \mathcal{E}_{T} ( r(X^\e,u^n)) \int_0^{T\wedge \tau^n} |u_s|_U^2 \,ds\Big)$ is non decreasing in $n$. Indeed, if $m<n$ we have $\tau_m \leq \tau_n$ almost surely, and therefore
\begin{align*}
& \E \Big( \mathcal{E}_{T} ( r(X^\e,u^n)) \int_0^{T\wedge \tau^n} |u_s|_U^2 \,ds\Big) 
=\E \Big( \mathcal{E}_{T\wedge \tau^n} (r(X^\e,u^n)) \int_0^{T\wedge \tau^n} |u_s|_U^2 \,ds\Big) \geq   \\ & \geq
\E \Big( \mathcal{E}_{T\wedge \tau^n} ( r(X^\e,u^n)) \int_0^{T\wedge \tau^m} |u_s|_U^2 \,ds\Big)
= 
\E \left(\E\Big( \mathcal{E}_{T\wedge \tau^n} ( r(X^\e,u^n)) \int_0^{T\wedge \tau^m} |u_s|_U^2 \,ds\Big | \mathcal{F}^W_{\tau_m\wedge T} \Big)\right )
 \\ &= 
\E \Big( \mathcal{E}_{T\wedge \tau^m} ( r(X^\e,u^n)) \int_0^{T\wedge \tau^m} |u_s|_U^2 \,ds \Big)=
 \E \left(\E\Big( \mathcal{E}_{T\wedge \tau^m} ( r(X^\e,u^m)) \int_0^{T\wedge \tau^m} |u_s|_U^2 \,ds\Big | \mathcal{F}^W_{\tau_m\wedge T} \Big)\right )
 \\ &= 
 \E \Big( \mathcal{E}_{T} ( r(X^\e,u^m)) \int_0^{T\wedge \tau^m} |u_s|_U^2 \,ds \Big).
 \end{align*}
 Thus, if $u\in \mathcal{U}^\e_{ad}$ then there exists finite  the limit  $\lim_{n\to +\infty} \E \Big( \mathcal{E}_{T} ( r(X^\e,u^n)) \int_0^{T\wedge \tau^n} |u_s|_U^2 \,ds\Big)  \in \mathbb{R}$.

\end{remark} 
\begin{remark}\label{remark-justification} To further justify the choice of the class $\mathcal{U}^\e_{ad}$ of admissible controls, recall that if we define 
$
d\,\mathbb{P}^{n} := \mathcal{E}_{T\wedge \tau^n} ( r(X^\e,u^n))d \mathbb{P}\;$ and 
$W^{n}_t := W_t-\int_0^t r(X^\e_s, u^n_s)ds,
$
then $(X^\epsilon)$ satisfies \eqref{stato eps u} with $(W)$ replaced by the $\mathbb{P}^{n}$ Wiener process $(W^{n}_t)_{t\geq 0}$. Namely:
 $$\left \{
     \begin{array}{ll}
        d X^{\e}_t =  A X^{\e}_t dt +b(X^{\e,u}_t) dt + \sigma( X^{\e}_t) (  \Gamma^\e [G (W^{n}+\int_0^{\cdot}r(X^{\e}_s,u_s)ds)] )_t dt, &  t \in (0,T], \\ \\
         X^\e(0)= x_0 ,& 
     \end{array} \right.$$
     and 
    $\displaystyle \E \Big( \mathcal{E}_{T\wedge \tau^n} ( r(X^\e,u^n)) \int_0^{T\wedge \tau^n} |u_s|_U^2 \,ds\Big)$ coincides with  $\displaystyle\E^{\mathbb{P}^{n}}  \Big( \int_0^{T\wedge \tau^n} |u_s|_U^2 \,ds \Big)$ and it seems natural to ask that :
    $\sup_n \displaystyle\E^{\mathbb{P}^{n}}  \Big( \int_0^{T\wedge \tau^n} |u_s|_U^2 \,ds \Big)<\infty.$
\end{remark}
\noindent We are eventually ready to rigorously introduce our cost functional $J^\e$. Namely, we set for any $u \in \mathcal{U}^\e_{ad}$:
\begin{equation}\label{funzionale costo}
    J^\e(x_0,u) :=  \lim _{n \to \infty} \E \Big[ \mathcal{E}_{T} ( r(X^\e,u^n)) \Big(\int_0^T l(X^\e_s,u^n _s)\, ds +  h(X^\e_T)\Big)\Big].
\end{equation}
Such functional is well defined, indeed we can prove:
\begin{proposition}\label{well posedness cost}
    For any $u \in \mathcal{U}^\e_{ad}$  the cost functional $J^\e(x_0, u)$  given in \eqref{funzionale costo} is well -defined.
\end{proposition}
\Dim 
We show that $ \E  \Big( \mathcal{E}_{T} ( r(X^\e,u^n))\int_0^T l(X^\e_s,u^n_s) \, ds \Big)$ 
is a (real-valued)  Cauchy sequence.
Let $m>n$. 
   By the martingale property of $ (\mathcal{E}_{t}( r(X^\e,u)))_t $  and \eqref{def g} we have: 
    \begin{align}\label{convergenza_costo}\begin{split}
        &\bigg|\E  \Big( \mathcal{E}_{T} ( r(X^\e,u^m))\int_0^T l(X^\e_s, u^m_s) \, ds \Big) 
        - 
        \E  \Big( \mathcal{E}_{T} ( r(X^\e,u^n))\int_0^T l(X^\e_s, u^n_s) \, ds \Big)\bigg|
        \\ & \leq 
        \bigg|\E \Big ( \,  \mathcal{E}_{T\wedge \tau_m} ( r(X^\e,u))\int_0^{ T \wedge \tau_n}l(X^\e_s,u^m_s) \, ds
        \Big)
        + 
        \E \Big ( \,  \mathcal{E}_{T\wedge \tau_m} ( r(X^\e,u))\int_{ T \wedge \tau_n}^Tl(X^\e_s,u^m_s) \, ds
        \Big) 
        \\ & \quad -
        \E \Big ( \,  \mathcal{E}_{T\wedge \tau_n} ( r(X^\e,u))\int_0^{ T \wedge \tau_n}l(X^\e_s,u_s) \, ds
        \Big)
        - 
        \E \Big ( \,  \mathcal{E}_{T\wedge \tau_n} ( r(X^\e,u))\int_{ T \wedge \tau_n}^Tl(X^\e_s,u_\star) \, ds
        \Big)\bigg|
         \\ & = 
         \bigg| \E \Big ( \,  \mathcal{E}_{T\wedge \tau_m} ( r(X^\e,u))\int_{ T \wedge \tau_n}^Tl(X^\e_s,u^m_s) \, ds
        \Big) 
        - 
        \E \Big ( \,  \mathcal{E}_{T\wedge \tau_n} ( r(X^\e,u))\int_{ T \wedge \tau_n}^Tl(X^\e_s,u_\star) \, ds
        \Big)\bigg|
        \\ & \leq 
        \bigg| \E \Big ( \,  \mathcal{E}_{T\wedge \tau_m} ( r(X^\e,u))\int_{ T \wedge \tau_n}^{ T \wedge \tau_m} l (X^\e_s,u_s) \, ds \Big)\bigg | 
        \\& \qquad 
        +
       T C_l(1+|u_\star|_U^2) 
        \left(
        \E (  \mathcal{E}_{T\wedge \tau_n} ( r(X^\e,u)) I_{\{\tau_n < T\}}) 
         +
        \E (  \mathcal{E}_{T\wedge \tau_m} ( r(X^\e,u)) I_{\{\tau_m < T\}})
        \right).
          \end{split} \end{align}
    We start from the the last two terms. It holds, by Markov inequality:

    \begin{align*}
    \E (  \mathcal{E}_{T\wedge \tau_n} ( r(X^\e,u)) I_{\{\tau_n < T\}})
       & =  \E (  \mathcal{E}_{T\wedge \tau_n} ( r(X^\e,u)) I_{ \{ \int_0 ^ {T \wedge \tau_n} |u_s|_U^ 2 \, ds \geq n\}}) \\
       & \leq  \frac{1}{n } \,   \E \Big(  \mathcal{E}_{T\wedge \tau_n} ( r(X^\e,u)) \int_0^{T\wedge \tau_n} |u_s|_U^2 \, ds  \Big),
\end{align*}       
        
and the same holds for $ \E (  \mathcal{E}_{T\wedge \tau_m} ( r(X^\e,u)) I_{\{\tau_m < T\}})$. In particular, since $u \in \mathcal{U}_{ad}$, both terms go to zero as $n,m \to \infty$.

 Regarding the first term,  we notice that it is smaller than 
and
 \begin{align*}
    &  \E \Big ( \,  \mathcal{E}_{T\wedge \tau_m} ( r(X^\e,u))\int^{ T \wedge \tau_m}_ { T \wedge \tau_n} |u_s|_U^2 \, ds \Big)
    \\ 
    &=  
    \E \Big ( \,  \mathcal{E}_{T\wedge \tau_m} ( r(X^\e,u))\int^{ T \wedge \tau_m}_ 0 |u_s|_U^2 \, ds \Big) 
    -  
    \E \Big ( \,  \mathcal{E}_{T\wedge \tau_n} ( r(X^\e,u))\int^{ T \wedge \tau_n}_ 0 |u_s|_U^2 \, ds \Big).
 \end{align*}
 Since $u \in \mathcal{U}_{ad}$ the sequence $ \Big (\E \Big ( \,  \mathcal{E}_{T\wedge \tau_n} ( r(X^\e,u))\int^{ T \wedge \tau_n}_ 0 |u_s|^2_U \, ds \Big) \Big)_n$ is a Cauchy sequence, see also Remark \ref{rem lim}. Therefore, the difference above as well converges to zero as $n,m \rightarrow 0$.

In a similar way we show that  $ \E  \Big( \mathcal{E}_{T} ( r(X^\e,u^n))h(X^\e_T) \Big)$ %
is a  Cauchy sequence.

\finedim

We can define the admissible controls and the cost functional of the limit control problem \eqref{stato u} in a similar way.

\section{BSDE Representation of the Value Function and of the Optimal Control}

\subsection{Hamiltonian function associated to the cost functional}\label{sec hamiltonian}
We introduce the Hamiltonian function $\psi$:
\begin{equation}
   \psi: K\times H^* \to \R, 
   \quad   
   \psi(x,z) := \inf _{u \in U} \{  l(x,u) - \langle z , r(x,u) \rangle \},
\end{equation}
where $\langle z , r(x,u) \rangle$ denotes the duality between $H$ and $H^*$.
Thanks to Assumptions \ref{hyp-control} we have the following:
\begin{corollary}\label{reg ham}
    The function $\psi$ has the following properties (for suitable constants $M_{\psi}$ and $L_{\psi}$):
    \begin{align}\label{prop psi}
   |\psi(x,z)|   &\leq  M_{\psi} (1+ |z|_{H^*}^2)  \qquad\qquad\qquad 
    \forall x \in K, \forall  z \in H^*, 
   \\
   |\psi(x,z)-\psi(x,z')| &\leq  L_{\psi}(1+ |z|_{H^*}+ |z'|_{H^*})  |z-z'|_{H^*} \qquad\qquad\,  \forall x \in K, \forall  z,z' \in H^*,
   \\
   |\psi(x',z)-\psi(x,z)| &\leq  L_{\psi}(1+ |z|_{H^*}^2) ( |x-x'|_K\wedge 1) \qquad\qquad\,  \forall x,x' \in K, \forall  z \in H^*.
    \end{align}
\end{corollary}

\Dim By \eqref{def g}  and \eqref{def r} we easily get that $\psi(x,z) \leq l(x,u_\star) - \langle z , r(x,u_\star) \rangle \leq  C_l |u_\star|_U^2$. 
On the other hand, there exists a finite constants $c$ such that, for every $u$ satisfying $|u|_U \geq c (1+ |z|_{H^*})$, it holds, recalling that $m_l>0$
\begin{align}
   l(x,u) - \langle z , r(x,u) \rangle 
   \geq 
   -c_l + m_l |u|_U^2 - M_r |z|_{H^*}(1 + |u|_U) 
  \geq 0,
\end{align}
while for $u$ satisfying $|u|_U \leq c (1+ |z|_{H^*})$ we have
\begin{align}
   l(x,u) - \langle z , r(x,u) \rangle 
   \geq 
   -c_l + m_l |u|_U^2 - M_r |z|_{H^*}(1 + |u|_U) 
   \geq -(c_l + cM_r + M_r)(1+|z|_{H^*}^2).
\end{align}
Hence, we deduce that there exists a costant $M_\psi$ such that 
\[  |\psi(x,z)|   \leq  M_\psi(1+ |z|_{H^*}^2) . \]
Next, the difference $|\psi(x,z)-\psi(x,z')| $ is controlled from above with
\begin{align*} 
    &\left| \inf_{|u|_U \leq c (1+ |z|_{H^*} +|z'|_{H^*})} (l(x,u) - \langle z , r(x,u) \rangle) 
    - 
    \inf_{|u|_U \leq c (1+ |z|_{H^*} +|z'|_{H^*})} (l(x,u) - \langle z' , r(x,u) \rangle) 
    \right| 
    \\ 
    &\qquad \leq \sup_{|u|_U \leq c (1+ |z|_{H^*} +|z'|_{H^*})} |z-z'|_{H^*} |r(x,u)| 
    \\
    &\qquad \leq  
    (c M_r + M_r) (1+ |z|_{H^*} +|z'|_{H^*}) |z-z'|_{H^*}
    \\
    &\qquad \leq  
    L_\psi (1+ |z|_{H^*} +|z'|_{H^*}) |z-z'|_{H^*}.
\end{align*}
Finally, in view of \eqref{def rr} and \eqref{def ll}, the difference $|\psi(x,z)-\psi(x',z)|$ is controlled from above with
\begin{align*}
      &\left| 
      \inf_{|u|_U \leq c (1+ |z|_{H^*})} (l(x,u) - \langle z , r(x,u) \rangle ) 
      - 
      \inf_{|u|_U \leq c (1+ |z|_{H^*})} (l(x',u) - \langle z , r(x',u) \rangle) \right|  
      \\  
      &\qquad  \leq \sup_{|u|_U \leq c (1+ |z|_{H^*})}  |l(x,u)-l(x',u)| 
      + 
      \sup_{|u|_U \leq c (1+ |z|_{H^*})} 
      |z|_{H^*} |r(x,u)-r(x',u)| 
      \\
      &\qquad \leq 
       2 c C_l|x-x'|_K (1+|z|_{H^*})^2
      + 
      c L_r|z|_{H^*}(1+|z|_{H^*})(	|x-y|_K\wedge 1)
      \\
      &\qquad \leq  L_{\psi}(1+ |z|_{H^*}^2) ( |x-x'|_K\wedge 1).
\end{align*}
\finedim

In the following, we assume that the infimum the definition of $\psi$ is indeed achieved.

\begin{Assumption}\label{argmin}
    There exists a measurable function $\underline{u}(x,z): K\times H^* \to U$ such that 
    \begin{enumerate}
        \item     
 $\psi(x,z)=\inf _{u \in U} \{  l(x,u) - \langle z , r(x,u) \rangle\}= 
        l(x,\underline{u}(x,z)) - \langle z , r(x,\underline{u}(x,z)) \rangle$;
        
        \item there exists a constant $ L_{\underline{u}} >0$, such that
        \begin{equation}\label{Lip control}
           | \underline{u}(x,z)-\underline{u}(x',z')|_U \leq L_{\underline{u}} \left[|z-z'|_{H^*}+(1+|z|_{H^*})(1\wedge |x-x'|_K)\right] \qquad \forall x, x'\in K, \forall z,z' \in H^* 
\end{equation}
    \end{enumerate}
\end{Assumption}
\begin{example} \label{esempio-feedback}
Assume that $U=H$ and let
$l(x,u):=l_0(x)+|u|_K^2$ and $r(x,u):=r_0(x)u$
with $l_0$ and $r_0$ bounded continuous functions $K\rightarrow \mathbb{R}$ and $K\rightarrow \mathcal{L}(H)$, respectively. 
In this case, if one identifies $H^*$ with $H$ by the canonical Riesz isomorphism, one gets
$$\psi(x,z)=l_0(x)-\frac{1}{4}|r_0(x)^*z|_{H^*}^2
\quad \hbox{and} \quad 
\underline{u}(x,z) = \frac12 r_0(x)^* z.
$$
Thus Assumptions \ref{argmin} are verified.
\end{example}

\subsection{BMO Martingales} \label{sec BMO} For the reader's convenience and in order to fix the notation, we report here a few basic facts on BMO martingales, following \cite{Kazamaki} and \cite{BriConf}.

Let $T \in (0,\infty)$ be given. A continuous  $(\Omega, (\mathcal{F}_t)_{t \in [0,T]},\mathbb{P})$ local martingale is a $BMO_2$ martingale on the time interval $[0,T]$ if 
$$ \Vert M \Vert_{BMO_2}
:= 
\sup_{\tau \in \mathcal{T}}\left\Vert \mathbb{E}\left( (M_T-M_{\tau})^2\big | \mathcal{F}_{\tau}\right)^{1/2}\right \Vert_{L^{\infty}(\Omega)}
= 
\sup_{\tau \in \mathcal{T}}\left\Vert \mathbb{E}\left( \langle M\rangle_T-\langle M\rangle_{\tau} \big | \mathcal{F}_{\tau}\right)^{1/2}\right \Vert_{L^{\infty}(\Omega)}
< \infty,$$
where $\tau$ in the supremum varies in the class $\mathcal{T}$ of all stopping times satisfying $\tau \leq T$ almost surely.
If $(\Psi)$ is a process in $L^{2,loc}_W(\Omega\times [0,T]; H^*)$ and $M_t=\int_0^t \Psi_s dW_s$, then 
\begin{equation}\label{M-BMO}
 \Vert M \Vert_{BMO_2}^2
 = 
 \sup_{\tau \in \mathcal{T}}\left\Vert \mathbb{E}\left( \int_{\tau}^T |\Psi_s |_{H^*}^2 ds \, \Bigg\rvert \, \mathcal{F}_{\tau}\right)\right \Vert_{L^{\infty}(\Omega)},
\end{equation} 
whenever the right-hand side is finite.

Moreover, again in the particular case $M_t=\int_0^t \Psi_s dW_s$, by \cite[p. 26]{Kazamaki} (see also \cite[Formula (13), p. 831]{BriConf}) one has that for all $p\geq 1$ there exists a finite constant $c(p)$ such that
\begin{equation}\label{stimaBMOLp}
    \mathbb{E}\left(\int_0^T |\Psi_s|_{H^*}^2 ds\right)^p\leq c(p) \Vert M \Vert_{BMO_2} ^{2p}.
\end{equation}
Finally, the exponential martingale 
$$\mathcal{E}(\Psi)_t:=\exp\left( \int_0^t \Psi_s dW_s - \frac{1}{2}  \int_0^t |\Psi_s|_{H^*}^2 ds  \right) $$
is  uniformly integrable and, by \cite[Formula (6), p. 824]{BriConf}, there exists $q^*>1$, depending only on $\Vert M \Vert_{BMO_2}$, such that for all $q\in (1,q^*)$ there is a suitable finite constant $C(q, \Vert M \Vert_{BMO_2})$ such that for every stopping time $\tau \leq T$ it holds.  
\begin{equation}\label{stimaBMOesponenziale}
    \mathbb{E}(\mathcal{E}(M)_T^q |\mathcal{F}_{\tau} )\leq C(q, \Vert M \Vert_{BMO_2})\mathcal{E}(M)_{\tau}^q.
\end{equation}
In particular, taking $\tau=0$ one gets
\begin{equation}\label{stimaBMOesponenzialesemplice}
    \mathbb{E}(\mathcal{E}(M)_T^q )\leq C(q, \Vert M \Vert_{BMO_2}).
\end{equation}

\subsection{BSDE representation}$ $
We are in the position to prove that:

\begin{theorem} \label{teo esist} 
Under Assumptions \ref{ass-state} and \ref{hyp-control}  here exists a unique triple of stochastic processes $(X^\e,Y^\e,Z^\e)$ {adapted to the filtration $(\mathcal{F}^W_t)_{t \in [0,T]}$} such that $X^\e $ has continuous trajectories, $Y^\e  \in  L^{\infty}_W(\Omega;C([0,T];\R)) $, $Z^\e \in  L^2_W(\Omega \times [0,T];H^*)) $, and $(X^\e,Y^\e,Z^\e)$ is a solution to the following system:
\begin{equation}\label{stato eps new} \left \{
     \begin{array}{ll}
        \frac{d}{dt} X^\e_t = ( A X^\e_t + b(X^\e_t))  + \sigma( X^\e_t)  (\Gamma ^\e[G W])_t  &  t\in (0,T], \\ \\
        -d Y^\e_t = \psi(X^\e_t, Z^\e_t) \, dt - Z^\e_t \, d W_t
        \\   \\ X^\e(0)= x_0, \quad Y^\e_T=h(X^\e_T).
     \end{array} \right.
 \end{equation}
Moreover
\begin{equation}\label{stimaBMO_eps}
 \sup_{t\in [0,T]} |Y^\e_t|_{ L^{\infty}_W(\Omega;C([0,T];\R))}+   \Big\Vert \int_0^{ \cdot}Z^\e_s dW_s  \Big \Vert _{BMO_2} \leq \kappa
\end{equation}
where $\kappa >0$ is independent of $\e$.
\end{theorem}
\Dim
By Theorem \ref{esistenza stato} the forward equation has a unique solution $X^\e \in L^2_W(\Omega;C([0,T];K))$.  Following \cite{BriConf}[Prop 11], see also \cite{Kobilanski} [Prop 2.1] there exists a  unique $(Y^\e,Z^\e)$, such that
{
\begin{equation}\label{stime stand Y Z eps}
   \sup_{t\in [0,T]} |Y^\e_t|_{ L^{\infty}_W(\Omega;C([0,T];\R))} \leq  M_{h} +  M_\psi T  
\end{equation}
\begin{equation}
     \E \int_0^T |Z^\e|^2_{H^*} \, dt  \leq C
\end{equation}}
for a constant $C>0$ depends only on $M_h, M_\psi, T$.

Let us check that \eqref{stimaBMO_eps}. 
We follow again \cite{BriConf}.

We apply  the It\^o formula to $\phi(Y_t+m)$, where $m$ is chosen so that $Y_t+m \geq 0$ and $\phi (x)=  (e^{2Cx} - 2Cx-1) / (2C^2)$, so that for all $x \geq 0$  satisfies $\phi'(x) \geq 0$ and $\frac{1}{2} \phi''(x)- C\phi'(x)=1$, for some $C > M_\psi$, given in \eqref{prop psi}. Then  taking the  conditional expectation with respect to $\F_\tau$, for any stopping time $\tau \leq T$, we have that- we avoid the subscript in the norms for simplicity:
\begin{equation*}
    \phi(Y_\tau +m) + \frac{1}{2}  \E^{\F_{\tau}} \Big( \int_\tau^T \phi''(Y_s+m)|Z_s|^2 \,ds  \Big)=  \E^{\F_{\tau}}  \phi(Y_T +m)+
    \E^{\F_{\tau}} \Big( \int_\tau^T \phi'(Y_s+m) \psi(Z_s)\,ds  \Big)
\end{equation*}
thus for every $\tau \leq T$:
\begin{equation*}
\begin{split}
     \phi(Y_\tau +m) + \frac{1}{2} \E^{\F_t} \Big( \int_\tau^T  |Z_s|^2 \,ds  \Big)& =  \E^{\F_t}  \phi(Y_T +m)+
    \E^{\F_t} \Big( \int_\tau^T  \phi'(Y_s+m) [\psi(Z_s)-C |Z_s|^2]\,ds  \Big)\\
    & \leq  \E^{\F_t}  \phi(Y_T +m)+ M_{\psi}
    \E^{\F_t} \Big( \int_\tau^T  \phi'(Y_s+m)\,ds  \Big)
\end{split}  
\end{equation*}
A similar argument is used to prove \eqref{stime stand Y Z eps}, see \cite{Kobilanski}[Prop 2.1].
And the claim holds by  \eqref{stime stand Y Z eps} and \eqref{M-BMO}.
\finedim

\begin{theorem}
  Under Assumptions \ref{ass-state}, \ref{hyp-control} and \ref{argmin} we have that:
   \begin{equation}\label{representation}
  Y^\e_0= \inf_{u \in \mathcal{U}_{ad}^\e} J^\e(u) = J^\e(\underline{u}(X^\e,Z^\e)).
\end{equation}  
where $(X^\e,Y^\e,Z^\e)$ is given in Theorem \ref{teo esist} and $\underline{u}$ is defined in Assumption \ref{argmin}.
\end{theorem}
\Dim
We need to get a fundamental relation for a generic $u \in \mathcal{U}^\e_{ad}$, or at least for its approximations. 

Proceeding as in Remark \ref{remark-justification} we apply Girsanov transformation to ensure that  $$W^n_t:= - \int_0^t r(X^\e_s,u^n_s) \, ds + W_t$$ is a cylindrical Wiener process under the probability $d\mathbb{P}^n:= \mathcal{E}_T(r(X^\e,u^n)) d \mathbb{P}$.
 Thus
\begin{align*}
   -d Y^\e_t & = \psi(X^\e_t, Z^\e_t) \, dt- Z^\e_t \, d W_t    =
    [\psi(X^\e_t, Z^\e_t) -Z^\e_s r(X^\e_t, u^n_t)] \, dt  - Z^\e_t \, d W^n_t 
\end{align*}
 Adding and subtracting the current cost and integrating between $0$ and $T$ we have:
\begin{align*}
    Y^\e_0 & = \E \Big( \mathcal{E}_T(r(X^\e,u^n)) \int_0^T [\psi(X^\e_t, Z^\e_t) -Z^\e_s r(X^\e_t, u^n_t) - l(X^\e_t, u^n_t)] \, dt \Big) +J^{\e,n}(u)
\end{align*}
where $  J^{\e,n}(u) : =  \E \Big[ \mathcal{E}_{T} ( r(X^\e,u^n)) \Big(\int_0^T (l(X^\e_s,u^n _s))\, ds +  h(X^\e_T)\Big)\Big]$.

The definition of $\psi$ yields $ J^n(u) \geq Y_0$ for every and consequently, by definition \eqref{funzionale costo},
\begin{equation} \label{stima costo facile}
    J^\e(u)= \lim_{n \to \infty} J^{\e,n}(u) \geq Y^{\e}_0 \qquad \qquad \forall u \in \mathcal{U}^\e_{ad}.
\end{equation}
Now we define  $ \underline{u}^\e(s)= \underline{u}(X^\e_s,Z^\e_s)$,  where $\underline{u} $ is given in Assumption \ref{argmin} and $(X^\e,Z^\e)$ is the solution of \eqref{stato eps new}.
From Assumption \ref{argmin} we have that 
\begin{equation}\label{stima contr ottimo}
    |\underline{u}^\e(s)| \leq C_{\underline{u}} (1+|Z^{\epsilon}_s|) 
\end{equation}    
for some constant $C_{\bar{u}}$.

We have to show that
 $\underline{u}^\e \in \mathcal{U}^\e_{ad}$ and 
     $Y^\e_0= J^e(\underline{u}^\e)$. To this purpose we define (see Definition \ref{Controlli Ammissibili eps}):
\begin{enumerate}
  \item  $ \bar{\tau}_n=\inf{ \{  0 \leq t \leq T : \int_0^t |\underline{u}^\e(s)|^2 \, ds \geq n \}}$,
  \item $\bar{u}^n(s)= \underline{u}^\e(s) I_{[0,\tau_n]}(s) + u^\star  I_{(\tau_n,T]}(s)$,
  \item $\bar{W}^n_t: -\int_0^t r(X^\e_s, \bar{u}^n(s))ds+W_t$.
\end{enumerate}
The backward component in \eqref{stato eps new}
can be rewritten as:
 \begin{equation}
      \begin{cases}
   - d Y^\e_t =     \left[\psi(X^\e_t,Z^\e_t)-  Z^\e_t r(X^\e_t, \bar{u}^n(t))\right] dt  -   Z^\e_t d\, W^n_t, \qquad t \in [0,T] \\
  \quad  \  Y^\e_T = h(X^\e_T).
\end{cases}
 \end{equation}
 Thus, recalling point (1) in Assumption \ref{argmin}
 \begin{equation}\label{rel fond n}
     Y^\e_0=\E(\mathcal{E}_T(r(X^\e,\bar{u}^n) h(X_T^\e)) + \E\Big( \mathcal{E}_T(r(X^\e,\bar{u}^n)\int_0^T l(X^\e_s, \bar{u}^n(s)) \, ds\Big) .
 \end{equation}

 By \eqref{def r} and \eqref{stima contr ottimo} we have, $\mathbb{P}-a.s.  $ 
 $$ |r(X^{\e}_s, \bar{u}^n(s))|\leq M_r(1+C_{\underline{u}})+M_r C_{\underline{u}} |Z^\e_s|$$
thus by \eqref{stimaBMO_eps} :
$$
    \sup_{\e>0,n\in \mathbb{N}}  \Big\Vert \int_0^{ \cdot}r^*(X^{\e}_s, \bar{u}^n(s)) dW_s  \Big \Vert _{BMO_2} <+\infty
$$
and finally the above estimate together with 
\eqref{stimaBMOesponenzialesemplice} yield that there exists $q>1$ such that:
\begin{equation}\label{stima q esp}
  \sup_{\e>0,n\in \mathbb{N}}  \E ( \mathcal{E}_T(r(X^\e,\bar{u}^n))^q) < \infty
\end{equation}

Again, by \eqref{stimaBMO_eps} and \eqref{stimaBMOLp}  we obtain that, for  every $p \geq 1$ :
\begin{equation}\label{stima Z p}
   \sup_{\e>0}   \E \left(\int_0^T |Z^\e_s|^2 \, ds \right)^{p/2} < \infty.
\end{equation}
Summing up we have that  $ \lim_{n \to +\infty} \bar{u}^n _s= \underline{u}^\e_s= \underline{u}(X^\e_s,Z^\e_s)$  for all $s \in [0,T],  \ \mathbb{P}-a.s.$.

Moreover for all $p\geq 1$
$$\begin{array}{rl}
    \E \left(\displaystyle\int_0^T | \bar{u}^n_s|^2 \, ds \right)^{p/2} &  \leq c(p)T^{p/2}|u^\star|^p+c(p) \E \left(\displaystyle\int_0^T | \underline{u}(X_s^\e,Z^\e_s)|^2 \, ds \right)^{p/2} \\ & \leq   \tilde{C} + \tilde{C}\E \left( \displaystyle\int_0^T |Z^\e_s|^2 \, ds \right)^{p/2}   < +\infty
\end{array}$$
 In view of \eqref{stima q esp} and of the above estimate the sequence $$\Big(\mathcal{E}_T(r(X^\e,\bar{u}^n))\displaystyle\int_0^T | \bar{u}^n_s|^2 \, ds\Big)_n$$
 turns out to be  uniformly integrable.  Moreover
$\int_0^T | \bar{u}^n_s|^2 \, ds\rightarrow \int_0^T | \underline{u}^\e_s|^2 \, ds $ $\mathbb{P}$-a.s. and $$ \mathcal{E}_T(r(X^\e,\bar{u}^n)) \rightarrow  \exp\left\{\int_0^T r(X^\e_s,\underline{u}^\e_s) dW_s - \frac{1}{2}\int_0^T |r(X^\e_s,\underline{u}^\e_s)|^2 ds\right\}
\; \hbox{$\mathbb{P}$-a.s.  }$$  thus the limit:
$\displaystyle \lim_{n \to \infty} \E\Big( \mathcal{E}_T(r(X^\e,\bar{u}^n) )\int_0^T|\bar{u}^\e_s|^2 \, ds   \Big)$ exists in $\mathbb{R}$
and we can conclude that $ \underline{u}^\e_s \in \mathcal{U}^\e_{ad}$.

In a similar way, taking into account  \eqref{def g} and \eqref{def h} we get that both $$\left(\mathcal{E}_T(r(X^\e,\bar{u}^n) h(X_T^\e))\right)_n \hbox{ and }  \left(\mathcal{E}_T(r(X^\e,\bar{u}^n)\int_0^T l(X^\e_s, \bar{u}^n(s)) \, ds\right)_n$$ 
are uniformly integrable and  $\mathbb{P}$-a.s. converging sequences of random variables. Letting $n\rightarrow \infty$ in  \eqref{rel fond n} we have that:
\begin{equation*}
        Y^\e_0=   \lim _{n \to \infty} \E \Big[ \mathcal{E}_{T} ( r(X^\e,\bar u^n)) \Big(\int_0^T l(X^\e_s,\bar u^n _s)\, ds +  h(X^\e_T)\Big)\Big]=J^\e(\underline{u}^\e)
\end{equation*}
and the claim follows. \finedim 

\medskip 
 Following exactly the same argument we prove that:
  \begin{theorem} \label{teo esist limite}
 { Under Assumptions  \ref{ass-state}, \ref{hyp-control} and \ref{argmin} we have here exists a unique solution $( \color{violet}{\hat{X}},\hat Y,\hat Z)$ with $ \hat{X}  \in  L^2_W(\Omega;C([0,T];K)), \hat Y  \in  L^2_W(\Omega;C([0,T];\R)) $, $\hat Z \in  L^2_W(\Omega \times [0,T];H^*)) $ to:
\begin{equation}\label{stato eps new lim} \left \{
     \begin{array}{ll}
     d \hat{X}_t = ( A \hat{X}_t + \hat{b}(\hat{X}_t)) \, dt + \sigma( \hat{X}_t) G  dW_t, &  t\in (0,T], \\ \\
        
        -d \hat Y_t = \psi(\hat X_t, \hat Z_t) \,  dt - \hat{Z}_t \, d W_t
        \\   \\   X(0)= x_0, \qquad  \hat Y_T=h(\hat X_T). &
     \end{array} \right.
 \end{equation}}
  Moreover
  \begin{equation}\label{stimaBMO}
 \sup_{t\in [0,T]} |\hat{Y}_t|_{ L^{\infty}_W(\Omega;C([0,T];\R))}+   \Big\Vert \int_0^{ \cdot}\hat{Z}_s dW_s  \Big \Vert _{BMO_2} \leq \kappa.
\end{equation}  Finally if
    \begin{equation}\label{def_hat_u}
    \hat{u}_t:= \underline{u}(\hat{X}_t, \hat{Z}_t)
    \end{equation}
     where  $\underline{u}$ is given in Assumption \ref{argmin} then, $\hat u$ is admissible (that is belongs to $\mathcal{U}_{ad}$) and 
\begin{equation*}
    \hat Y_0 = \inf_{u \in \mathcal{U}_{ad}}\hat J(u) = \hat J(\hat{u}).
\end{equation*}
where $$ \hat J(u) : = \lim_{n\rightarrow \infty}  \E \Big[ \mathcal{E}_{T} ( r(\hat X, \hat u^n)) \Big( h(\hat X_T)+\int_0^T (l(\hat X_s,\hat u^n _s))\, ds  \Big)\Big]$$ and 
  $ \hat{\tau}_n=\inf{ \{  0 \leq t \leq T : \int_0^t |\hat{u}(s)|^2 \, ds \geq n \}}$, $\hat{u}^n(s)= \hat{u}(s) I_{[0,\hat \tau_n]}(s) + u^\star I_{(\hat \tau_n,T]}(s)$.
\end{theorem}

\section{Limit problem and convergence}\label{sec-conv}
We have developed all the machinery necessary to approach the convergence of the control problems. 
Of course, this convergence is only expected when the maps $(\Gamma^\e)$ are ``good'' approximations of the noise. 
However, it turns out that only convergence of the forward equation alone is necessary, and no information about the control problems has to be assumed.
More precisely, let us assume the following natural condition:
{\begin{Assumption}\label{ipotesi limite forward} Let $X^\e$ and $X$ the solutions to \eqref{stato eps new} and \eqref{stato eps new lim} respectively.
We assume that for every $t\in [0,T]$,  $X^\e_t \to \hat X_t $ in probability as $\e \to 0$.
 \end{Assumption}}

We are now in a position to prove our main convergence result.
\begin{theorem}\label{main_genenral}
    Under Assumptions  \ref{ass-state}, \ref{hyp-control}, \ref{argmin} and \ref{ipotesi limite forward}, we have that:
    \begin{equation}\label{limite valore}
      \lim_{\e \to 0}\inf_{u \in \mathcal{U}_{ad}}J^\e(u)=\inf_{u \in \mathcal{U}_{ad}}\hat J(u)
    \end{equation}
    Moreover if  $\underline{u}^\e$ and  $\hat{u}$ are the optimal  admissible controls introduced in the previous section then as we know $J^\e(\underline{u}^\e)=\inf_{u \in \mathcal{U}_{ad}}J^\e(u)$;
     $\hat{J}(\hat{u})=\inf_{u \in \mathcal{U}_{ad}}\hat{J}(u)$, moreover:
    \begin{equation}\label{limite controlli ottimi}
         \lim_{\e \to 0}  \E \int_0^T |\underline{u}^\e_s-\hat {u}_s|^2\, ds =0.
    \end{equation}
\end{theorem}
\Dim 
Let us consider the equation for the difference $Y^\e_t-\hat{Y}_t:=\tilde{Y}^\e_t$. It solves :
\begin{align}\label{eq diff}
    \tilde{Y}^\e_t= h(X^\e_T)-  h(\hat{X}_T)+\int_t^T (\psi(X^\e_s,Z^\e_s)-\psi(\hat{X}_s,\hat{Z}_s))\, ds + \int_t^T 	\tilde{Z}^\e_s d W_s
\end{align}
where $ \tilde{Z}^\e_t=  Z^\e_t-Z_t$. Equation \eqref{eq diff} can be rewritten as:
\begin{align}
    \tilde{Y}^\e_t= h(X^\e_T)-  h(\hat{X}_T) + \int_t^T (\psi (X^\e_s, \hat Z_s)- \psi (\hat X_s, \hat Z_s)) \, ds + \int_t^T  K^\e_s \tilde{Z}^\e_s\, ds + \int_t^T\tilde{Z}^\e_s\ d W_s
\end{align}
where \begin{equation}
    K^\e_s=: v(X^\e_s, Z^\e_s,\hat Z_s) \hbox{ and }  v(x,z,z') = \left\{  \begin{array}{ll}
        \frac{\psi(x,z)-\psi(x,z')}{|z-z'|^2}(z-z')&  \hbox{ if }|z-z'| \not=0 \\
       0  &   \hbox{ if }|z-z'|=0 
    \end{array}\right.
\end{equation}
Notice that  $ K^\e_s\leq L_{\psi}(  1+ |Z^\e_s|+ |\hat Z_s|)$ thus in view of
 by  \eqref{stimaBMO_eps} and \eqref{stimaBMO}
 $$\sup_{\e>0}\left\Vert \int_0^\cdot K^\e_s dW_s\right\Vert_{BMO_2}<+\infty$$ 
 Moreover if $f^\e_s:= \psi (X^\e_s, \hat Z_s)- \psi (\hat X_s, \hat Z_s)$ then in view of \eqref{prop psi}
 \begin{equation}\label{stima-f}
     |f^\e_s|\leq L_{\psi}(1+|\hat Z_s|^2)(1\wedge |X^\e_s-\hat{X}_s|).
 \end{equation}
By \eqref{stimaBMOLp} and \eqref{stimaBMO} we have that
\begin{equation}\label{stimaf} \sup_{\e >0}\mathbb{E}\left(\int_0^T |f^\e_s|ds\right)^q <+\infty, \hbox{ for all $q\geq 1$}
\end{equation} 
 thus assumption A3 in \cite{BriConf} is verified for any $p>1$. Consequently we can apply estimate (7) in \cite{BriConf} with $p^*=2p$

\begin{equation}\label{stimafinale}
\begin{split}
 & ( \E \sup_{t \in [0,T]}|\tilde{Y}^\e_t|^p)^{1/p} + \Big( \E( \int_0^T|\tilde{Z}^\e_t|^2\, dt )^{p/2} \Big)^{1/p} \\& \leq  C\Big \{[ \E(|h(X^\e_T)- h(\hat{X}_T)| )^{2p}] ^{1/2p} + \Big[\E \Big( \int_0^T |f^\e_s| \, ds  \Big)^{2p}\Big]^{1/2p}\Big\}\Big(1+ \Big[\E\Big(\int_0^T |K^\e_s|^2\, ds \Big)^{3p/2}\Big]^{1/3p}\Big) \leq \\& \leq  \tilde{C}\Big \{[ \E(|h(X^\e_T)- h(\hat{X}_T)| )^{2p}] ^{1/2p} + \Big[\E \Big( \int_0^T |f^\e_s| \, ds  \Big)^{2p}\Big]^{1/2p}\Big\}
\end{split}
\end{equation}
where $\tilde{C}$ depends on $||\int_0^{\cdot}K^\e_s dW_s||_{BMO_2}$ see again \cite{BriConf}. 

Moreover, recalling that $\int_0^T |\hat{Z}_s|^2ds\in L^q$ for all $q\geq 1$ we readily deduce that the sequence $\Big( \int_0^T |f^\e_s| \, ds  \Big)^{2p}$ is uniformly integrable. To prove that $\mathbb{E}\left(\int_0^T |f^\e_s|ds\right)^{2p} \rightarrow 0 $   it is therefore enough to prove that $\int_0^T |f^\e_s|ds$ converges to 0 in probability.

We start by showing that, for almost every $s\in [0,T]$, $\mathbb{E} |f^\e_s| \rightarrow 0$. By \eqref{stima-f} it is enough to show that $\mathbb{E}(1+|\hat Z_s|^2)(1\wedge |X^\e_s-\hat{X}_s|) \rightarrow 0$. Indeed, in view of Assumption \ref{ipotesi limite forward}, $(1+|\hat Z_s|^2)(1\wedge |X^\e_s-\hat{X}_s|))$ converges to 0 in probability and is dominated by $(1+|\hat Z_s|^2)$.

Then again by dominated convergence $\mathbb{E} \int_0^T |f^\e_s| ds \rightarrow 0$ and consequently  $\int_0^T |f^\e_s| ds \rightarrow 0$  in probability.

In the same way taking into account assumption (Hp 2-h) we get that $\E(|h(X^\e_T)- h(\hat{X}_T)| )^{2p}\rightarrow 0$ thus by \eqref{stimafinale}
\begin{equation}\label{limite Y e Z}
    \lim_{\e \to 0} \E \sup_{t \in [0,T]}|\tilde{Y}^\e_t|^p=0 \qquad  \lim_{\e \to 0}  \E\left( \int_0^T|\tilde{Z}^\e_t|^2\, dt \right)^{p/2}  =0
\end{equation}
and we deduce \eqref{limite valore}.

 It remains  to prove  \eqref{limite controlli ottimi}. Recalling that $\underline{u}^\e_s=\underline{u}(X^\e_s,Z^\e_s)$ and  $\hat{u}_s=\hat{u}(\hat{X}_s,\hat{Z}_s)$ by \eqref{Lip control} we have:
 $$\mathbb{E}\int_0^T |\underline{u}^\e_s-\hat{u}_s|^2ds\leq L_{\underline{u}}^2\mathbb{E}\int_0^T|Z^\e_s-\hat Z_s|^2ds +  L_{\underline{u}}^2 \mathbb{E}\int_0^T(1+|\hat Z_s|^2)(1\wedge  |X^\e_s-\hat X_s|)^2ds
$$ 
Then \eqref{limite controlli ottimi} follows by Assumption \eqref{ipotesi limite forward} and relation \eqref{limite Y e Z} exactly as in the above detailed  proof that  $\mathbb{E}\left(\int_0^T |f^\e_s|ds\right)^{2p} \rightarrow 0 $.
\finedim

\begin{remark}
By \cite[Theorem 2.2]{AsFlaPap}, our motivating example \eqref{statopap}-\eqref{statopapreduced} satisfies Assumptions  \ref{ass-state} and \ref{ipotesi limite forward} when $K$ is a finite dimensional Hilbert space. 
Therefore, our \autoref{main_genenral} applies as soon as the cost functional satisfies Assumptions \ref{hyp-control} and \ref{argmin}, and we have the convergence of the optimal costs and controls associated control problem \eqref{statopapcontr} towards those of \eqref{statopapreducedcontr}.
\end{remark}


\section{Examples and further developments} \label{sec:Examples}

\subsection{Wong-Zakai type approximations} $ $

Let us consider a simple finite-dimensional stochastic equation
 \begin{equation} \label{WZstato} \left \{
     \begin{array}{ll}
        d_t X_t =  \sigma( X^\e_t) dB_t,   &  t \in [0,T], \\ 
         X(0)= x_0,
     \end{array} \right.
 \end{equation}
where $(B_t)_{t\in [0,T]}$ is a $\mathbb{R}^d$-valued standard Brownian motion,  $\sigma: \mathbb{R}^n \rightarrow L(\mathbb{R}^{ d},\mathbb{R}^{n})$ is a regular map.
We define the usual It\^o-Stratonovich correction map $\sigma_2: \mathbb{R}^n \rightarrow L(\mathbb{R}^d\times \mathbb{R}^d,\mathbb{R}^n) $ by
$$\sigma_2(x)(v,w)=\nabla_x(\sigma (x)v)(\sigma (x)w),\qquad \hbox{for all}\; x\in \mathbb{R}^n, v,w\in \mathbb{R}^d  $$
Following \cite{CaGuZa} we assume that $\sigma $ and $\sigma_2$ are of class $C^3$ with bounded first second and third derivative.

\medskip

Concerning regularization of noise let $\rho: \mathbb{R}\rightarrow [0, +\infty)$ be a smooth function $\rho$ with compact support  satisfiying $\rho(s)=0, \hbox{ for $s<0$ and } \int_{-\infty}^0\rho(s) ds =1.$ For all $\e>0$ let $\rho_\e(s):=\e^{-1}\rho(\e^{-1} s)$.

Given a continuous semimartingale $I$ in the class $\mathcal{I}$ introduced in  paragraph \ref{settings} let 
$$ I^{\e}_t =(\rho^\e \star I)_t\quad\hbox{and}\quad \Gamma^{\e}[I]_t =\dot{I^{\e}_t}\quad\hbox{in particular}\quad B^{\e}_t =(\rho^\e \star I)_t\quad\hbox{and}\quad \Gamma^{\e}[B]_t =\dot{B^{\e}_t}$$
where $I$ has be extended to $0$ before $0$ and after $T$. Notice that Assumption \ref{ipotesi su Gamma} is satisfied due to the asymmetry of the mollifier $\rho$. 

Concerning the control problem let 
$(u_t)_{t\in [0,T]}\in L^2_B([0,T];\mathbb{R}^m)$ and $r$, $l$, $h$ satisfy Assumption \ref{hyp-control} with $K=\mathbb{R}^n  $ and $U=\mathbb{R}^m$.
We consider the approximating controlled equations 
 \begin{equation} \label{wz-contr1e} \left \{
     \begin{array}{ll}
        d_t X^{u,\e}_t =  \sigma(X^{u,\e}_t)\Gamma^\e[ \int_0^{\cdot}r(X_s,u_s) ds  + B]_t dt,   &  t \in [0,T], \\ 
         X(0)= x_0,
     \end{array} \right.
 \end{equation}
 which can be rewritten as:
  \begin{equation*} 
  \left \{
     \begin{array}{ll}
        d_t X^{u,\e}_t =  \sigma(X^{u,\e}_t)
        [\rho_e \star r(X,u)]_t dt +  \sigma(X^{u,\e}_t)\dot{B}^{\e}_t dt,   &  t \in [0,T], \\ 
         X(0)= x_0,
     \end{array} \right.
 \end{equation*}
 and the limit controlled equation
  \begin{equation*} 
  \left \{
     \begin{array}{ll}
        d_t X^{u}_t = \textrm{Tr}[\sigma_2 (X^{u}_t)]+ \sigma(X^{u}_t) r(X^{u}_t,u_t) dt + \sigma(X^{u}_t)d{B}_t,   &  t \in [0,T], \\ 
         X(0)= x_0,
     \end{array} \right.
 \end{equation*}
 together with the cost functionals $J^e$ and $J$ defined as in \eqref{def funzionale controllo}.
 In \cite{CaGuZa} it is shown that, under the present assumptions, if $X^\e$ solves
  \begin{equation} \label{wz-e}\left \{
     \begin{array}{ll}
        d_t X^{\e}_t =  \sigma(X^{\e}_t)\dot{B}^{\e}_t dt,   &  t \in [0,T], \\ 
         X(0)= x_0,
     \end{array} \right.
 \end{equation} 
 and $\hat{X}$ solves
  \begin{equation} \label{wz-hat}\left \{
     \begin{array}{ll}
        d_t \hat{X}_t = \textrm{Tr}[\sigma_2 (X^{u}_t)]dt+ \sigma( \hat{X}_t)dB_t,   &  t \in [0,T], \\ 
         X(0)= x_0,
     \end{array} \right.
 \end{equation} 
 then $X^{\e}\rightarrow \hat{X}$, $\mathbb{P}-a.s.$ in a suitable Holder norm and in particular   $X^{\e}_t\rightarrow \hat{X}_t$, $\mathbb{P}-a.s.$ for all $t>0$.

 Thus if Assumption \ref{argmin} holds we are in a condition to apply Theorem \ref{main_genenral} and conclude that 
  \begin{equation*}
      \lim_{\e \to 0}\inf_{u \in L^2_B([0,T];\mathbb{R}^m)}J^\e(u)=\inf_{u \in L^2_B([0,T];\mathbb{R}^m)}\hat J(u)
    \end{equation*}
    Moreover  there exist optimal controls   $\underline{u}^\e$ and  $\hat{u}$ in $L^2_B([0,T];\mathbb{R}^m)$ such that $$J^\e(\underline{u}^\e)=\inf_{u \in L^2_B([0,T];\mathbb{R}^m)}J^\e(u)\quad \hbox{and}\quad \hat{J}(\hat{u})=\inf_{u \in L^2_B([0,T];\mathbb{R}^m)}\hat{J}(u).$$
    Finally  $\displaystyle \lim_{\e \to 0}  \E \int_0^T |\underline{u}^\e_s-\hat {u}_s|^2\, ds =0.$

\subsection{Quadratic fast-fast interaction}
A possible extension of our results could take into account climatic systems with fast-fast interaction at meteorologic scales, replacing \eqref{statopap} with
 \begin{equation} \label{statopap_q} \left \{
     \begin{array}{ll}
        d X^\e_t = A X^\e_t dt + {b}(X^\e_t) dt  + \sigma( X^\e_t) Q^\e_t dt,   &  t \in [0,T], \\ 
         X^\e(0)= x_0, & \\
         d Q^\e_t= q(Q^\e_t,Q^\e_t) dt -\dfrac{1}{\e}  Q^\e_t dt + \dfrac{1}{\e} Gd{W}_t, &  t \in [0,T], \\ 
         Q^\e(0)=0. &
     \end{array} \right.
 \end{equation}
In the equation above, $q:H \times H \to H$ is a continuous bilinear map.
For simplicity we suppose that $K$ is finite dimensional.
Hereafter we shall implicitly assume conditions on $q$ guaranteeing existence and uniqueness of solutions to \eqref{statopap_q} for a sufficient class of noises $W$.

Technically speaking, the results in \cite{MTV,AsFlaPap} do not cover the case of quadratic self-interaction for the fast variable and therefore require $q=0$. In view of their geophysical interpretation, assuming $q=0$ is a restrictive modelling assumption (cf. Equation 2.4 in \cite{MTV}) since most equations of geophysical fluid dynamics do have quadratic nonlinearities.
These difficulties have been recently overcome in a series of papers \cite{FlPa21,FlPa22,DePa}.

A stochastic model reduction of \eqref{statopap_q} is performed in \cite{DePa}, where convergence in probability towards a reduced equation is proved. The reduced equation has the form 
 \begin{equation} 
 \left \{
     \begin{array}{ll}
        d\hat{X}_t = A \hat{X}_t dt +  \hat{b}(\hat{X}_t) dt  + \sigma( \hat{X}_t)GdW_t + \sigma(\hat{X}_t) \hat{q} dt,  &  t \in [0,T], \\ 
         \hat{X}(0)= x_0,
     \end{array} \right.
 \end{equation}
 where $\hat{q}$ is the average of the fast-fast interaction with respect to the centered Gaussian measure with covariance $Q:=\frac12 G^* G$, namely
 \begin{align*}
\hat{q} := \int_H q(w,w) \mathcal{N}(0,Q)(d w).
 \end{align*}

In view of the results of this paper, one can introduce the controlled fast-slow system with quadratic fast-fast interaction
 \begin{equation} \label{statopapcontr_q} \left \{
     \begin{array}{ll}
        d X^{\e,u}_t = A X^{\e,u}_t dt+  {b}(X^{\e,u}_t) dt  + \sigma( X^{\e,u}_t) Q^{\e,u}_t dt,   &  t \in [0,T], \\ 
         X(0)= x_0, & \\
         d Q^{\e,u}_t= q(Q^{\e,u}_t,Q^{\e,u}_t) dt -\dfrac{1}{\e}  Q^{\e,u}_t dt +\dfrac{1}{\e} Gr(u_t)dt+  \dfrac{1}{\e} Gd{W}_t, &  t \in [0,T], \\ 
         Q(0)=0, &
     \end{array} \right.
 \end{equation}
 and the controlled reduced equation
 \begin{equation} \label{statopapreducedcontr_q} \left \{
     \begin{array}{ll}
        d\hat{X}^u_t = A \hat{X}^u_t dt+  \hat{b}(\hat{X}_t) dt + \sigma( \hat{X}^u_t)Gr(u_t)dt  + \sigma( \hat{X}^u_t)GdW_t  + \sigma( \hat{X}^u_t)\hat{q} dt,  &  t \in [0,T], \\ 
         \hat{X}(0)= x_0.
     \end{array} \right.
 \end{equation}

Convergence of the optimal control problems falls into our general theory by considering the maps
\begin{align}
    \Gamma^\e (I):= Q,
\end{align}

where $Q$ is the unique solution of 
\begin{align*}
    dQ_t = q(Q_t,Q_t) dt - \dfrac{1}{\e}  Q_t dt + \dfrac{1}{\e} dI_t.
\end{align*}

The analogue of \autoref{main_genenral} follows.

\end{document}